%% file: heat_tree_optimality.tex
\documentclass[a4paper,12pt]{article}
\usepackage[a4paper,margin=1in,headheight=\baselineskip]{geometry}
\usepackage[utf8]{inputenc}
\usepackage[T1]{fontenc}
\usepackage{microtype,comment}
\usepackage[english]{babel}
\usepackage[hidelinks,hypertexnames=false]{hyperref}
\usepackage{graphicx}
\usepackage{amsmath}
\usepackage{amssymb}
\usepackage{amsthm}
\usepackage{cleveref}
\usepackage{thmtools}
\usepackage{mathtools}
\usepackage{enumerate}
\usepackage{enumitem}
\usepackage{color}
\usepackage{esint}
\usepackage{algorithm}
\usepackage{algpseudocode}
\usepackage[section]{placeins}
\usepackage{stmaryrd} 
\usepackage{subcaption}
\usepackage{bbm}
\usepackage{tikz}
\usetikzlibrary{arrows.meta}

\numberwithin{equation}{section}

\usepackage{scrlayer-scrpage}
\ofoot{}

\crefname{subsection}{subsection}{subsections}
\crefname{subsubsection}{subsubsection}{subsubsections}

\input{newcommands}

\def\tend{{t_{\rm end}}}

\title{Instance optimality of a space-time adaptive Radau method}
\author{Michael Feischl\textsuperscript{1}, Benedikt Gr\"a\ss le\textsuperscript{2}, David Niederkofler\textsuperscript{1}\footnote{The main ideas and proof arguments are human generated and we used AI for the technical details and implementation. We provide a partial LEAN certificate of the results at \url{https://github.com/michaelFeischl/space-time-optimality-LEANv2}. The certificate proves the $n=2$ case of Theorem~\ref{gen:instance} under the assumptions that certain external references are correct. Funded by the Deutsche Forschungsgemeinschaft (DFG, German Research Foundation) -- Project-ID 258734477 -- SFB 1173, the Austrian Science Fund (FWF) under the special research program Taming complexity in PDE systems (grant SFB F65) as well as projects I6667-N and P36150. This work was also supported by the ERC grant OPTIMAL (grant agreement No. 101125225). }\\[0.5em]
\small\textsuperscript{1}Institute of Analysis and Scientific Computing, TU Wien\\
\small Wiedner Hauptstraße 8--10, 1040 Vienna, Austria\\[0.3em]
\small\textsuperscript{2}Institut für Mathematik, Universität Zürich\\
\small Winterthurerstr.~190, CH-8057 Zürich, Switzerland}

\hypersetup{pdftitle={Instance optimality of a space-time adaptive Radau method}}
\begin{document}
\raggedbottom
\maketitle
\begin{abstract}
We prove instance optimality for the residual estimator of a space--time
Radau finite element method for the heat equation. The algorithm first uses D\"orfler marking to reduce
the estimator by a constant factor. It then constructs a new time partition and new conforming
spatial meshes by coarsening of the fine Galerkin solution.
\end{abstract}

\section{Introduction}\label{sec:introduction}
The theory of optimal adaptive mesh refinement originated from breakthrough results by Binev, Dahmen, DeVore~\cite{bdd}, Stevenson~\cite{stevenson07}, and Cascon, Kreuzer, Nochetto, Siebert~\cite{ckns}, who showed that a standard adaptive loop of the form
\begin{align*}
    \fbox{Solve}\longrightarrow\fbox{Estimate}\longrightarrow\fbox{Mark}\longrightarrow\fbox{Refine}
\end{align*}
produces optimal convergence rates for the Poisson problem. 
The new ideas inspired a flurry of research in this area, extending the results to many other model problems, see e.g.,~\cite{ks,cn} for conforming methods,~\cite{BeMao10,ncstokes3}    
for nonconforming methods,~\cite{LCMHJX,CR2012} for mixed formulations, and~\cite{fkmp,gantumur} for boundary integral equations.
For a comprehensive overview, we refer to~\cite{carstensen_2014,actaadaptive}. 

One of the earliest works~\cite{bdd} used an additional coarsening step to obtain the optimality results, i.e.,
\begin{align*}
    \fbox{Solve}\longrightarrow\fbox{Estimate}\longrightarrow\fbox{Mark}\longrightarrow\fbox{Refine}\longrightarrow\fbox{Coarsen}
\end{align*} This coarsening step was later used to show optimality of $hp$-adaptive FEM in~\cite{hpopt}.
The authors argue that this coarsening step is necessary for $hp$-refinement, as one can construct examples were pure refinement misses optimal meshes. The basic idea is
to follow the standard adaptive loop until one measures reduction of the estimator by a fixed factor. Then, one uses the coarsening routine to trim over refinement from the mesh.

We follow this route for space-time adaptivity. Our recent work~\cite{FHN2026} forms the basis for this algorithm, as it already shows linear convergence of an adaptive time-refinement loop for Radau timestepping for the heat equation. We extend these results to obtain linear convergence of adaptive space-time refinement for time slab meshes, i.e., adaptive time steps with possibly different adaptive spatial meshes attached to each time step. We use the inf-sup stability of the underlying discretization together with new quasi-orthogonality results~\cite{qogeo,infsup} to show that the standard adaptive loop for adaptive time-slab meshes leads to a reduction of the error estimator after a fixed number of iterations. Then, we show that a coarsening algorithm, that reapproximates the current discrete solution starting from the initial mesh based on~\cite{inst:BD} finds instance optimal space-time meshes in almost linear cost.

Our method is restricted to adaptive time-slab meshes with a certain temporal grading condition, i.e., the time steps must not increase too quickly. While fully unstructured adaptive meshes would be attractive from an approximation theory point of view, time slabs have the big advantage that the method can still be used in classical time stepping fashion without a large memory overhead.

Other methods for space-time adaptivity with theoretical guarantees use non-standard discretizations of the time-dependent problem, e.g.,~\cite{diening2025} where the authors find a symmetric reformulation of the heat equation in non-standard Sobolev spaces, or~\cite{FUHRER202127,GaSte24} where a least squares reformulation is used for adaptivity.  
We also mention~\cite{adaptiveG,adaptiveDG}, which develop adaptive time stepping algorithms for parabolic problems that are guaranteed to reach the final time and obey a given error tolerance.

Adaptive wavelet discretizations~\cite{wavelet_schwab} are known to be optimal for quite some time now but come with challenging implementation and other restrictions. See also~\cite{halfwavelet} for a hybrid, wavelet-in-time approach.

Section~\ref{inst:section} defines the discretization, estimator and
algorithm. Section~\ref{sec:optimality} states the main result and
describes the proof. Section~\ref{sec:proof-overview} gives an overview of the proof. Section~\ref{gen:numerics} compares adaptive and
uniform refinement for the heat equation on an L-shaped domain.

\section{Model problem and algorithm}\label{inst:section}

\subsection{Heat equation and admissible meshes}
Let $\Omega\subset\R^n$, $n\ge1$, be a bounded polyhedral Lipschitz domain,
$\tend>0$, $H=L^2(\Omega)$ and $V=H^1_0(\Omega)$ with
$\|v\|_V=\|\nabla v\|_{L^2(\Omega)}$. We consider
\begin{equation}\label{model:heat}
 \partial_tu-\Delta u=f\quad\text{in }(0,\tend)\times\Omega,
 \qquad u=0\quad\text{on }(0,\tend)\times\partial\Omega,
 \qquad u(0)=u_0 .
\end{equation}
We suppress the time interval in Bochner space notation. The weak
solution $u\in\X=H^1(V^*)\cap L^2(V)$ has initial trace $u_0\in H$
and satisfies
$\langle\partial_tu,v\rangle+(\nabla u,\nabla v)=(f,v)$ for all
$v\in V$ at almost every time. We use
$\|u\|_{\X}^2=\|\partial_tu\|_{L^2(V^*)}^2+
\|u\|_{L^2(V)}^2+\|u(0)\|_H^2$.
Here $V^*=H^{-1}(\Omega)$ is the dual of $V$ in its gradient norm.
We write $\PP q(T)$ for spatial polynomials of total degree at most
$q$, and $\PP q(I;Z)$ for degree-at-most-$q$ temporal polynomials
with coefficients in a Hilbert space $Z$. The symbols $\lesssim$
and $\simeq$ allow constants depending on the fixed data of the
discretization, but not on refinement levels or local mesh sizes.

Fix polynomial degrees $p\ge1$, $k\ge2$ and set $m=\max\{2,6-k\}$ (for technical reasons, we will need to $m$-sect every time interval
on refinement).
The initial mesh $\QQ_0$ consists of a partition $\I_0$ of $(0,\tend]$
and one spatial triangulation per time interval. All its spatial
triangulations belong to a common simplicial bisection forest
generated from a fixed initial spatial mesh~\cite{inst:Stevenson,inst:BFV}.
Each temporal subdivision replaces one interval by $m$
equal children. The children inherit the parent's spatial mesh,
which may be refined independently on each child.

For a time partition $\I_\QQ$ and conforming spatial meshes $\T_I$,
write
\[
 \QQ=\{I\times T:I\in\I_\QQ,\ T\in\T_I\},\qquad
 \#\QQ=\sum_{I\in\I_\QQ}\#\T_I .
\]
At the right endpoint of $I$, set $\T_I^+=\T_I\LCC\T_{I^+}$,
where $I^+$ is the successor and $\LCC$ denotes largest common
coarsening. On the last slab put $\T_I^+=\T_I$. We write
$\T\vee \T'$ for the coarsest common refinement of two spatial meshes.
The order $\T\preceq \T'$ means that $\T'$ refines $T$; we use the
same convention for space--time meshes and temporal partitions.
For fixed $\Const{grad}\ge1$ and $\Const{rCFL}>0$ require
\begin{equation}\label{model:admissible}
  |I^+|\le\Const{grad}|I|\quad(\text{if }I^+\text{ exists}),\qquad
  \max_{T\in\T_I^+}\operatorname{diam}(T)^2\le\Const{rCFL}|I|
       \quad(I\in\I_\QQ).
\end{equation}
We assume that the initial mesh $\QQ_0$ satisfies these inequalities. Write $\mathbb Q(\QQ_0)$
for its admissible refinements by the stated subdivisions with
conforming spatial meshes and \eqref{model:admissible}. 
For any admissible $\QQ$, define $\mathbb Q(\QQ)$ in the same way
with $\QQ$ as the starting mesh.
Given an error estimator $\eta$ (see next section), we measure the best attainable estimator at a prism budget $N$ by
\begin{equation}\label{inst:best}
 \sigma(N)=\inf\{\eta(\QQ):\QQ\in\mathbb Q(\QQ_0),\ \#\QQ\le N\},
 \qquad \inf\varnothing=+\infty .
\end{equation}

\subsection{Radau solution and error estimators}
Let $S(\T)=\SS p_0(\T)=V\cap\{v:v|_T\in\PP p(T),\ T\in \T\}$ and
let $P_\T$ be the spatial $L^2$ projection onto $S(\T)$. With
$h_\T|_T=|T|^{1/n}$ we assume the uniform estimates
\begin{equation}\label{model:projectors}
 \|h_\T^{-1}P_\T v\|_H\le C_{\rm stb}\|h_\T^{-1}v\|_H\quad(v\in H),
 \qquad \|P_\T v\|_V\le C_{\rm stb}\|v\|_V\quad(v\in V)
\end{equation}
all admissible spatial meshes. We note that for a wide range of polynomial degrees and spatial dimensions, this has been proved
in~\cite{DieGrading1,DieGrading2}.

Assume $f\in L^2(H)$ and $f|_I\in H^1(I;V^*)$ for every
$I\in\I_0$. The optimality theorem below assumes $u_0\in S(\T_{I_1^0})$,
where $I_1^0$ is the first initial slab.

On $I=(a,b]$ let $r_1<\cdots<r_k=b$ denote the right Radau nodes,
$R_I$ their degree-$(k-1)$ interpolation operator, and $\Pi_I^j$
the temporal $L^2$ projection onto $\PP j(I)$. The discrete solution
$u_\QQ$ is continuous in time, satisfies
$u_\QQ|_I\in\PP k(I;S(\T_I))$ and takes the value
$P_{\T_{I_1}}u_0$ at zero, with $I_1$ the first slab of $\QQ$.
Its value at $b$ belongs to $S(\T_I^+)$. On each slab solve
\[
 (\partial_tu_\QQ(r_i),v)_H+(\nabla u_\QQ(r_i),\nabla v)
       =((\Pi_I^k f)(r_i),v)_H
 \quad\begin{cases}
    v\in S(\T_I),&1\le i<k,\\
    v\in S(\T_I^+),&i=k .
 \end{cases}
\]
The right endpoint value supplies the incoming value on the next slab.

Write $\|\cdot\|_D=\|\cdot\|_{L^2(D)}$. Let $\mathcal E_I$ be the
$L^2(I)$ projection onto the orthogonal complement of
$\{q\in\PP k(I):q(b)=0\}$ in $L^2(I)$; it retains both the
polynomial endpoint component and the nonpolynomial remainder.
For $Q=I\times T$, let $T^+\in\T_I^+$
contain $T$ and put $h=|T|^{1/n}$, $H_Q=|T^+|^{1/n}$.
Set $r=f-\partial_tu_\QQ+\Delta u_\QQ$ elementwise, and on an
interior face $L$ set $j_L=[\nabla u_\QQ]\cdot\nu_L$, with either
choice of the unit normal. Define the indicators by
\begin{align}
 d(I)&=|I|^2\|(1-\Pi_I^{k-2})\partial_t f\|_{L^2(I;V^*)}^2,
 \nonumber\\
 \eta_t(I)^2&=\|(1-R_I)\nabla u_\QQ\|_{I\times\Omega}^2+d(I),
 \label{inst:time-ind}\\
 \eta_x(Q)^2&=|T|^{2/n}\|r\|_Q^2+(|T^+|^{2/n}-|T|^{2/n})\|\mathcal E_Ir\|_Q^2
 \nonumber\\
 &\quad+\sum_{L\subset\partial T\cap\Omega}
  \big(|T|^{1/n}\|j_L\|_{I\times L}^2+(|T^+|^{1/n}-|T|^{1/n})\|\mathcal E_Ij_L\|_{I\times L}^2\big).
 \label{inst:space-ind}
\end{align}
The sum runs over interior faces of $T$. The global estimator is
$\eta(\QQ)^2=\sum_{I\in\I_\QQ}\eta_t(I)^2+
\sum_{Q\in\QQ}\eta_x(Q)^2$. Collective indicators always denote
the square root of the sum of their squared local contributions;
in particular $\eta_x(\QQ)$ and $\eta_t(\I_\QQ)$ have this meaning.
Face integrals use $(n-1)$-dimensional Hausdorff measure.

\subsection{Space-time refinement}
Given $\QQ$ and $0<\theta\le1$, pool the temporal and spatial
squared indicators. Mark in decreasing order until the selected
intervals $\mathcal M_t$ and prisms $\mathcal M_x$ satisfy
\begin{equation}\label{inst:bulk}
 \sum_{I\in\mathcal M_t}\eta_t(I)^2+
 \sum_{Q\in\mathcal M_x}\eta_x(Q)^2\ge\theta\eta(\QQ)^2 .
\end{equation}
Any choice satisfying this bulk criterion is sufficient.
The following refinement realizes the required strict reductions.
\begin{enumerate}
\item Whenever $I\in\mathcal M_t$ has an unmarked successor with
$m|I^+|>\Const{grad}|I|$, add that successor to $\mathcal M_t$ and
repeat to closure. Plan one $m$-section on each marked interval.
\item For each marked $I\times T$, bisect $T$ on $I$. Also bisect
its old reduced element $T^+$ wherever it is a leaf on $I$ or
$I^+$. Use conforming bisection closure for all these requests.
On the last slab only $I$ is involved.
\item Carry the resulting spatial meshes onto the proposed time
children. Enforce the second inequality of \eqref{model:admissible}.
\end{enumerate}
On every descendant of a marked prism, both its spatial cell and
its reduced endpoint cell are strict descendants of $T$ and $T^+$,
respectively, on their region. These two strict refinements enter
the reduction estimate. Refining a time slab alone leaves its
inherited spatial elements available for later spatial refinement.

\subsection{Space-time coarsening}
The work~\cite[Algorithm~1]{inst:BFV} propose a variant of the original Binev tree coarsening algorithm that results in conforming meshes.
Their algorithm stops after a given number of bisections. For the present coarsening in space and time, however, we use
finite weights from the fine Radau solution and load, stop all local
searches on different time slabs at a common price, and choose the time partition by a shortest
path argument in a graph of possible time intervals (see below) that penalizes mismatched traces.
Figure~\ref{inst:fig-coarsening} illustrates the procedure described in detail below.

\input{coarsen}

\paragraph{Spatial weights.}
Let $\F$ be the fine mesh after nested refinement, with solution $u_{\F}$.
For $s>0$, let $\T^{\mathrm{floor}}(s)$ be the coarsest uniform refinement of the initial spatial mesh
with ${\rm diam}(T)^2\leq \Const{rCFL}s$ for all $T\in\T^{\mathrm{floor}}(s)$.

For each interval $I$ of its temporal refinement history of $\F$, start the search at
$\T_I^{\mathrm{floor}}=\T^{\mathrm{floor}}(|I|/\Const{grad})$.
Let $D_T$ and $\pi_T^p$ be the pointwise-in-time $L^2(T)$ projectors
onto $\nabla\PP p(T)$ and $\PP p(T)$, respectively. Set
\begin{equation}\label{gen:sp-weight}
 e_I(T)=\|(1-D_T)\nabla u_{\F}\|_{I\times T}^2
       +|T|^{2/n}\|(1-\pi_T^p)f\|_{I\times T}^2,
 \qquad E_I(\T)=\sum_{T\in \T}e_I(T).
\end{equation}
The search sums only contributions for which $T$ was bisected in
forming the spatial mesh $\T_J^{\F}$ on fine slab $J$:
\begin{equation}\label{gen:finite-weight}
 \widehat e_I(T)=\!
 \sum_{\substack{J\in\I_{\F},\ J\subseteq I\\
                   T\text{ internal node in the refinement tree of }\T_J^{\F}}}
 e_J(T),\qquad
 \widehat E_I(\T)=\sum_{T\in \T}\widehat e_I(T).
\end{equation}

\paragraph{Spatial search with budget $B$ and price $\beta$.}
For each coarsening trial in the main algorithm (Algorithm~\ref{inst:algorithm}) below, fix $B$ and $\beta$. Perform the spatial search below for each $I$ in the temporal refinement history of $\F$.

For an element $T$, the priority with local weight $\widehat e_I(T)$ proposed in~\cite{inst:BFV} is
\begin{equation}\label{inst:priority}
 \mathfrak p(T)=\big(\nu(T)+1/\widehat e_I(T)\big)^{-1}
 \quad\text{if }\widehat e_I(T)>0,\qquad \mathfrak p(T)=0\quad\text{otherwise}.
\end{equation}
Following~\cite{inst:BFV}, we call a patch of spatial elements an \emph{executable necessary patch}, if it can be bisected conformingly and must
be bisected in every conforming closure that bisects $T$. For example, such a patch exists at the end of every mesh closure refinement chain. 

\begin{algorithmic}[1]
\State Skip $I$ for this budget if $\#\T_I^{\mathrm{floor}}>B$.
\State Initialize the mesh at $\T_I^{\mathrm{sel}}=\T_I^{\mathrm{floor}}$ with
       $W=\#\T_I^{\mathrm{floor}}$ and $\nu(\cdot)=0$, or resume its saved state in case the search already ran with a smaller budget $B$.
\While{$\max\mathfrak p>\beta$}
  \State If $W+1>B$, save the state and skip $I$ for this budget.
  \State Mark a leaf $T$ of maximum priority and choose one executable
         necessary patch.
  \State Increase $\nu(T)$ by $1/\widehat e_I(T)$.
  \State Bisect the patch, let each child $T''$ of a bisected parent $T'$ inherit $\nu(T''):=\nu(T')$,
         and increment $W$.
\EndWhile
\State Save the state and retain $I$ with the corresponding selected mesh $\T_I^{\mathrm{sel}}$,
       cost $W_I=W$, and leaf sum $\widehat E_I(\T_I^{\mathrm{sel}})$.
\end{algorithmic}

With at most $\Lambda$ cells per patch, we have $W\le\#\T\le\Lambda W$.

\paragraph{Time partition and realization.}
We consider the temporal error on each interval
and the trace mismatch at each interface $t=\sup I=\inf J$:
\begin{align}
 L_{\F}(I)&=\|(1-\Pi_I^{k-1})\nabla u_{\F}\|_{I\times\Omega}^2+d(I),
 \nonumber\\
 K_{\F}(I,J)&=\min(|I|,|J|)\|\Pi_I^k u_{\F}(t)-\Pi_J^k u_{\F}(t)\|_V^2,
 \nonumber\\
 K_{\F}(0,I)&=|I|\|u_0-\Pi_I^k u_{\F}(0)\|_V^2\quad(\inf I=0).
 \label{gen:practical-weights}
\end{align}
For a partition $\mathcal I=(I_1,\ldots,I_r)$ by retained intervals, the
score is
\begin{align}
 T_{\F}(\mathcal I)&=K_{\F}(0,I_1)+\sum_{i=1}^r L_{\F}(I_i)
                    +\sum_{i=1}^{r-1}K_{\F}(I_i,I_{i+1}),\nonumber\\
 \Phi_\beta(\mathcal I)&=T_{\F}(\mathcal I)
          +\sum_{i=1}^r\big(\widehat E_{I_i}(\T_{I_i}^{\mathrm{sel}})+\beta W_{I_i}\big).
 \label{gen:price-path}
\end{align}
We build a graph consiting of all temporal intervals in the refinement tree of $\F$.
To each interval $I$, assign the node cost $L_{\F}(I)+\widehat E_I(\T_I^{\mathrm{sel}})+\beta W_I$.
We introduce edges $I\to J$ of cost $K_{\F}(I,J)$ to consecutive admissible intervals, i.e., when $\sup I=\inf J$ and
$|J|\le\Const{grad}|I|$. We add source edges of cost $K_{\F}(0,I)$ at time
zero and sink edges of cost zero at $\tend$. 
Counting the node and edge costs, we search a cheapest path from $t=0$ to $t=\tend$ by increasing the right endpoint.
 with cost $+\infty$ if none exists (this can happen if the corresponding spatial meshes are too expensive).

On the chosen path write $\T_i^{\mathrm{sel}}=\T_{I_i}^{\mathrm{sel}}$
and repeat the first and last cores outside the path. Overlay consecutive
cores so both slabs at each interface contain the same two cores:
\begin{equation}\label{gen:sp-overlay}
 \T_i^{\mathcal U}=\T_{i-1}^{\mathrm{sel}}\vee \T_i^{\mathrm{sel}}\vee \T_{i+1}^{\mathrm{sel}},\qquad
 \mathcal U=\{I_i\times T:1\le i\le r,\ T\in \T_i^{\mathcal U}\}.
\end{equation}

\paragraph{Cost.}
For $N=\#\F$, coarsening takes $O(N\log^2(2+N))$ operations across
all budgets in the unit-cost tree model, reusing the spatial searches.
Count separately $O(N)$ queries each for $e_J(T)$ and $d(I)$, and the
solves and estimator evaluations. The score test bounds each trial
mesh by $3\Lambda B$ prisms, hence $O(N)$ prisms in total over the
doubled budgets.

\begin{algorithm}[ht]
\caption{Adaptive Radau method with conforming Binev coarsening}
\label{inst:algorithm}
\begin{algorithmic}[1]
\Require Admissible $\QQ_0$, bulk parameter $0<\theta\le1$ and
reduction factors $0<\lambda<\mu<1$.
\State Solve on $\QQ_0$ and compute $\eta_0=\eta(\QQ_0)$.
\For{$\ell=0,1,\dots$ while $\eta_\ell>0$}
 \State Set $\F=\QQ_\ell$. While $\eta(\F)>\lambda\eta_\ell$,
 mark by \eqref{inst:bulk}, apply nested refinement, solve and
 recompute $\eta(\F)$.
 \State Set $\tau=\mu\eta_\ell$ and $B=\#\QQ_0$. Prepare the time
 forest, spatial weights, and temporal costs from $\F$.
 \Repeat
  \State Set $\beta=\tau^2/B$. Run the spatial search with this price
  and cap $B$. Minimize \eqref{gen:price-path} on the retained intervals by shortest path search.
  \State If the minimum is at most $\tau^2$, form
  $\mathcal U$ by \eqref{gen:sp-overlay}, solve and compute $\eta(\mathcal U)$.
  \State If such a mesh satisfies $\eta(\mathcal U)\le\tau$, accept
  $\QQ_{\ell+1}=\mathcal U$ and set $\eta_{\ell+1}=\eta(\mathcal U)$.
  Otherwise double $B$.
 \Until{a mesh has been accepted}
\EndFor
\end{algorithmic}
\end{algorithm}

\section{Optimality theorem}\label{sec:optimality}
We state the main theorem of the work under the assumptions made above.
\begin{theorem}[Instance optimality in space and time]\label{gen:instance}
Let $\QQ_0$ be an admissible initial mesh of $(0,\tend]\times\Omega$
and $I_1$ its first time slab. Suppose that
\[
 f\in L^2(0,\tend;L^2(\Omega)),\qquad
 f|_I\in H^1(I;H^{-1}(\Omega))\quad(I\in\I_{\QQ_0}),\qquad
 u_0\in\SS p_0(\T_{I_1}).
\]
Fix $p\ge1$, $k\ge2$, $m=\max\{2,6-k\}$ and the one-sided time
grading and reverse CFL constants of Section~\ref{inst:section}.
Und the previous assumptions and for $0<\theta\le1$, $0<\mu<1$ there is a constant $K>1$ such
that the following conclusions hold for $0<\lambda\le\mu/K$.
Apply Algorithm~\ref{inst:algorithm} with bulk fraction $\theta$,
fine reduction factor $\lambda$ and measured acceptance factor $\mu$.
Every phase starting with $\eta_\ell=\eta(\QQ_\ell)>0$ terminates,
and its accepted mesh satisfies $\eta_{\ell+1}\le\mu\eta_\ell$.
There is a constant $C\ge1$ such that at every accepted level,
including the initial level,
\begin{equation}\label{gen:instance-eq}
 \boxed{\displaystyle
 \eta(\QQ_\ell)\le C
   \inf_{\substack{\QQ_{\mathrm{cmp}}\in\mathbb Q(\QQ_0)\\ C\#\QQ_{\mathrm{cmp}}\le\#\QQ_\ell}}
        \eta(\QQ_{\mathrm{cmp}}).}
\end{equation}
Here $\#\QQ_{\mathrm{cmp}}=\sum_{I\in\I_{\QQ_{\mathrm{cmp}}}}\#\T_I$ is the prism count.
By \eqref{inst:best}, the right-hand side equals
$C\sigma(\lfloor\#\QQ_\ell/C\rfloor)$ with $\inf\varnothing=+\infty$.
The constants may depend on the fixed initial discretization, the
domain, the polynomial degrees, the refinement family, the projection
stability and the fixed admissibility and algorithm parameters. They
are independent of the data, the refinement level and the number of
slabs.
\end{theorem}

\subsection{Overview of the proof}\label{sec:proof-overview}
The algorithm follows broadly the familiar refine and coarsen strategy first proposed in
~\cite{inst:BDD} for the Poisson equation.

We use  the new ideas for quasi-orthogonality in time dependent problems~\cite{inst:Feischl} to show that a fixed number of refinements reduce the error estimator by a given factor. Then, we use a modification of the conforming version of Binev's tree coarsening~\cite{inst:BFV} to coarsen the space-time mesh and to make sure
its size is minimal up to a constant factor.

The growing use of AI in mathematical research underscores the importance of transparent verification; in this spirit, we provide a partial Lean certificate for Theorem~\ref{gen:instance}, with its scope and assumptions stated explicitly.
We outline the main proof arguments below and reserve the full technical proof for a forthcoming paper.
A full Lean certificate was beyond the computational budget available to the authors, however, we believe that the
two dimensional case contains all major difficulties and we have high confidence in the few external references the LEAN certificate uses as assumptions.
See \url{https://github.com/michaelFeischl/space-time-optimality-LEANv2} for details on the certificate.
\paragraph{Standard adaptive loop, restarted after each coarsening.}
Fix an accepted mesh $\QQ_\ell$ and write
$\F_0=\QQ_\ell\preceq \F_1\preceq\cdots$ for the temporary refinement
sequence, $v_j=u_{\F_j}$.
Spatial marking refines both a marked cell and its reduced endpoint
cell; refining only the former would leave part of the endpoint
indicator unchanged. The resulting spatial contraction factor is
$2^{-1/(2n)}$, while temporal subdivision gives $m^{-1}$.
D\"orfler marking and the perturbation estimates yield
\[
 \eta_{j+1}^2\le a\eta_j^2+
 B\|v_{j+1}-v_j\|_\X^2,\qquad 0<a<1.
\]
We use the following precise form of relaxed quasi-orthogonality from~\cite{infsup,qogeo}:
\begin{equation}\label{overview:qo}
 \sum_{j=s}^{s+L-1}\|v_{j+1}-v_j\|_\X^2
 \le C_{\rm qo}(L-1)\eta_s^2,\qquad
 C_{\rm qo}(L-1)\lesssim L^{1-\delta},\quad\delta>0.
\end{equation}
The constant may grow with the block length, but it grows sublinearly.
Summing the estimator recurrence therefore gives
\[
 \min_{1\le i\le L}\eta_{s+i}^2
 \le \frac{a+B C_{\rm qo}(L-1)}{(1-a)L}\eta_s^2.
\]
A fixed number $L_\lambda$ of steps thus suffices to reach
$\eta(\F)\le\lambda\eta(\QQ_\ell)$. The argument requires no
minimal-cardinality marking and works for every $0<\theta\le1$.

\paragraph{Space-time modifications behind quasi-orthogonality.}
The main step is to show uniform inf-sup stability of the discretization. This is done in~\cite{FHN2026} for a time-adaptive Radau method with fixed spatial meshes and can be extended to the present setting.

Two features prevent a direct application of a fixed-form Galerkin
argument. The collocation test space changes its temporal shape when
a slab is split, and the discrete load uses the mesh-dependent
projection $\Pi_\mathcal I^k f$. We
construct a bounded adjoint lifting $R_\mathcal I^\dagger$ whose
range $D^\dagger(\QQ)$ is nested under refinement. Its functions
have degree $k-2$ on the $m$ children of each slab; the choice
$m=\max\{2,6-k\}$ supplies the moment and endpoint conditions
needed for this construction. The resulting formulation uses the
continuous parabolic bilinear form.

To obtain a fixed right-hand side as well, we introduce the
projected load as a second unknown in a nested polynomial space.
The resulting product problem is uniformly inf--sup stable, so the
finite-chain quasi-orthogonality result of~\cite{inst:Feischl} applies.
Its load increments are controlled by the temporal data term $d(I)$.
The common lifting $\chi(t)=(1-t/t_1)_+u_0$, supported on the first
initial slab $(0,t_1]$, reduces the trace to zero. Both the lifting
and its polynomial load contribution cancel in increments. Thus
$u_0$ only needs to be resolved on that first initial slab.
These steps give \eqref{overview:qo} in the full parabolic graph norm.

\paragraph{The coarsening step.}
The coarsening strategy follows~\cite{inst:BD} and its conforming variant~\cite{inst:BFV} for the spatial approximation but additionally needs to find an optimal time partition.  We need to make sure that changing meshes through time does not lead to large derivatives. Therefore, we include the time-interface terms in the coarsening search. 
Spatial coarsening considers the projections $\Pi_I^ku_{\F}$, which need not be continuous over time interfaces. 
By controling the temporal jumps of those projections, we can also control the quality of the continuous temporal approximation and thus obtain a quasi-optimal coarsend space-time mesh. Once we can find quasi-optimal spatial meshes for each allowed time interval, finding the optimal time partition is done by checking each admissible time interval starting from $t=0$ and working towards $\tend$.

\paragraph{Proof of instance optimality.}
This, again, is rather standard.
A competitor space-time mesh $\QQ_{\mathrm{cmp}}$ can refine time more deeply than $\F$ in some regions
and less deeply in others. The precise choice of constants in the coarsening algorithm ensures that if
$\#\QQ_{\mathrm{cmp}}\le B/K$ and $\eta(\QQ_{\mathrm{cmp}})\le\tau/K$ and
$\eta(\F)\le\tau/K$, the competitor (or a similar good one) is found.

\section{Numerical experiments}\label{gen:numerics}
The reference implementation of these experiments with scripts to recreate the plots can be found at \url{https://github.com/michaelFeischl/optimality_final}.
Let $\Omega=(-1,1)^2\setminus([0,1)\times(-1,0])$ be the L-shape and $\tend=1$.
We set $f=1$ and $u_0=1$ and impose homogeneous Dirichlet data on
$(0,1]\times\partial\Omega$.
The initial value belongs to $L^2(\Omega)$. On a mesh $\QQ$,
write $\T$ for the triangulation on its first slab and $P_\T$ for the
spatial $L^2$ projection onto $\SS 1_0(\T)$. We prescribe
$u_{\QQ}(0)=P_\T1$ by the consistent mass matrix. We use degree $p=1$
in space and $k=2$ in time. We measure the augmented error estimator 
\begin{equation}\label{gen:eta-aug}
 \eta_{\rm aug}(\QQ)^2
 = \sum_{I\in\I_{\QQ}}\eta_t(I)^2
   +\sum_{Q\in\QQ}\eta_x(Q)^2+\|1-P_\T1\|_{L^2(\Omega)}^2
\end{equation}
including the initial data oscillation term. This experiment with
$u_0=1$ is outside the scope of Theorem~\ref{gen:instance}.

The common initial mesh has one time slab and 24 triangles, obtained
by two uniform newest-vertex bisection generations of six root
triangles. Time refinement uses $m=4$ equal children, as prescribed
by $k=2$. We fix
\[
 \theta=0.5,\qquad \lambda=0.65,\qquad \mu=0.82,\qquad
 \Const{grad}=4,\qquad \Const{rCFL}=128.
\]
Thus consecutive slab lengths may increase by at most a factor of
four, and the reverse CFL condition in \eqref{model:admissible}
uses the constant 128. In the coarsening search the floor
uses $|I|/4$ and is never coarser than the initial triangulation.
The bulk fraction applies to the squared estimator, whereas
$\lambda$ and $\mu$ multiply the estimator itself. The latter
factors were chosen empirically, without estimating the constant
$K$ in Theorem~\ref{gen:instance}. In each phase, nested refinement
continues until $\eta_{\rm aug}(\F)\le0.65\eta_{\rm aug}(\QQ_\ell)$.
The budget search starts at $B=24$, sets
$\tau=0.82\eta_{\rm aug}(\QQ_\ell)$ and $\beta=\tau^2/B$, and
doubles $B$ after an unsuccessful trial. Each candidate passing
the score test is solved anew and accepted only if its recomputed
augmented estimator is at most $\tau$.

The marking is designed to minimize
$\#\mathcal M_t+\#\mathcal M_x$. It does not assert minimality of the
resulting prism count.

For the uniform comparison meshes define $\QQ_{g,j}$ by using $4^j$ uniform time slabs and $6\cdot2^g$
triangles per slab, $j\ge0$, $g\ge2$. Thus at budget
$N$ the uniform value is
\[
 E_{\rm unif}(N)=\min\{\eta_{\rm aug}(\QQ_{g,j}):
       6\cdot 2^g4^j\le N,\quad 2^{1-g}\le128\cdot4^{-j}\}.
\]

\begin{figure}
\centering
\includegraphics[width=\linewidth]{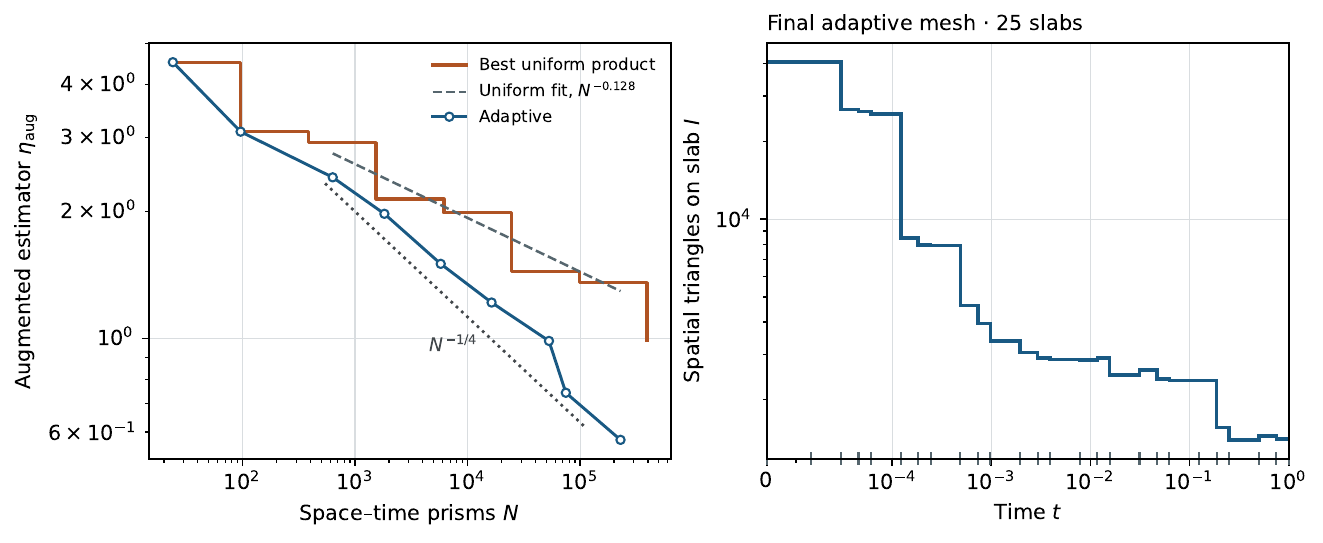}
\caption{$u_0=1$: The left panel shows
the augmented Radau estimator for the accepted meshes and the best
affordable uniform meshes. The right panel gives the triangle count
on each of the $25$ final slabs. Marks on the time axis denote slab
endpoints.}
\label{gen:fig-constant}
\end{figure}

\begin{figure}
\centering
\includegraphics[width=\linewidth,trim={0 1.7cm 0 0},clip]{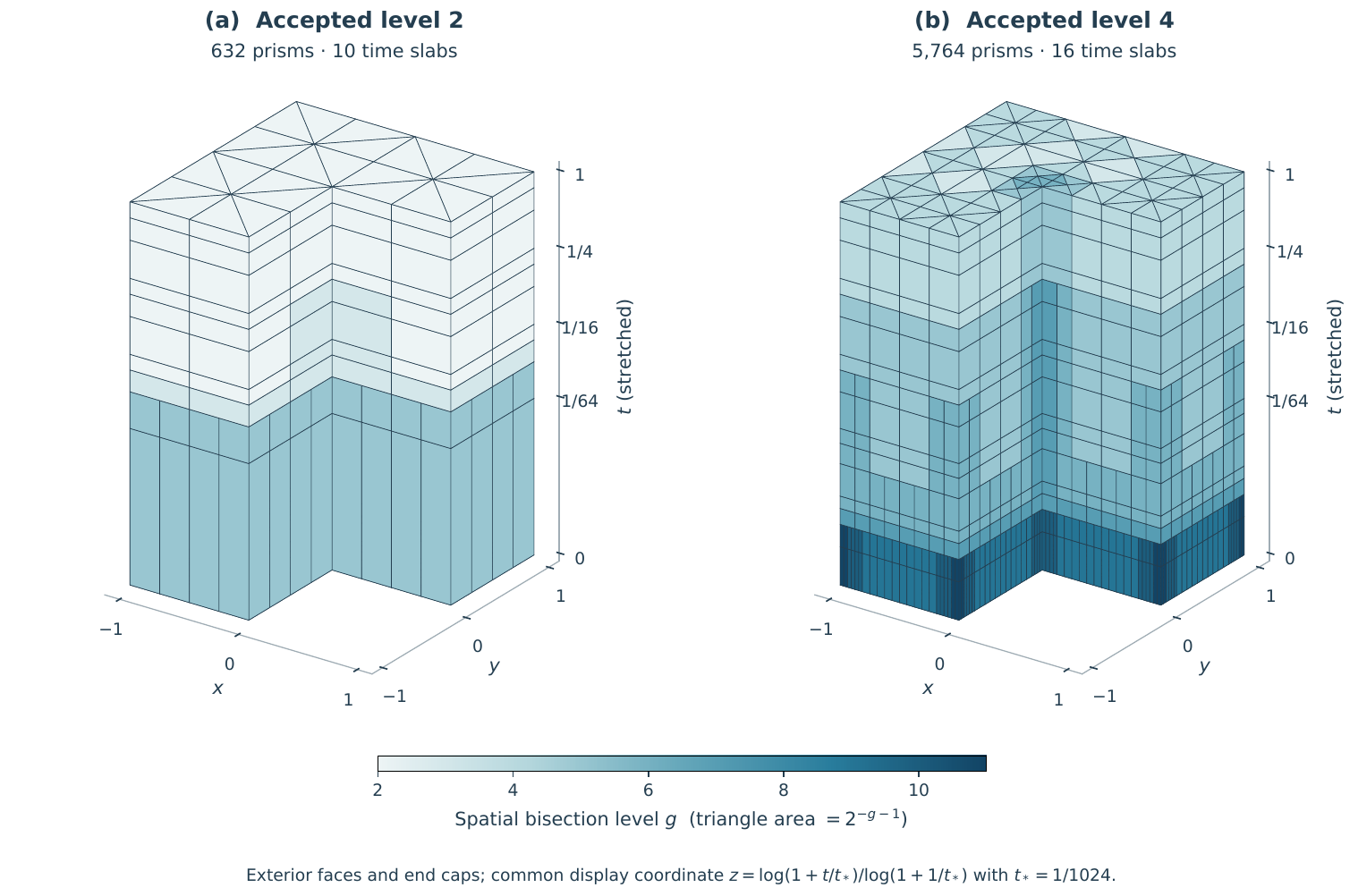}
\caption{$u_0=1$: Three-dimensional views of two accepted space--time meshes
for the constant-data example. The exterior prism faces and end caps
are shown; interior faces are hidden. Color encodes
the spatial bisection level, darker cells are finer.}
\label{gen:fig-spacetime-meshes}
\end{figure}

\begin{figure}
\centering
\includegraphics[width=\linewidth,trim={0 1.7cm 0 0},clip]{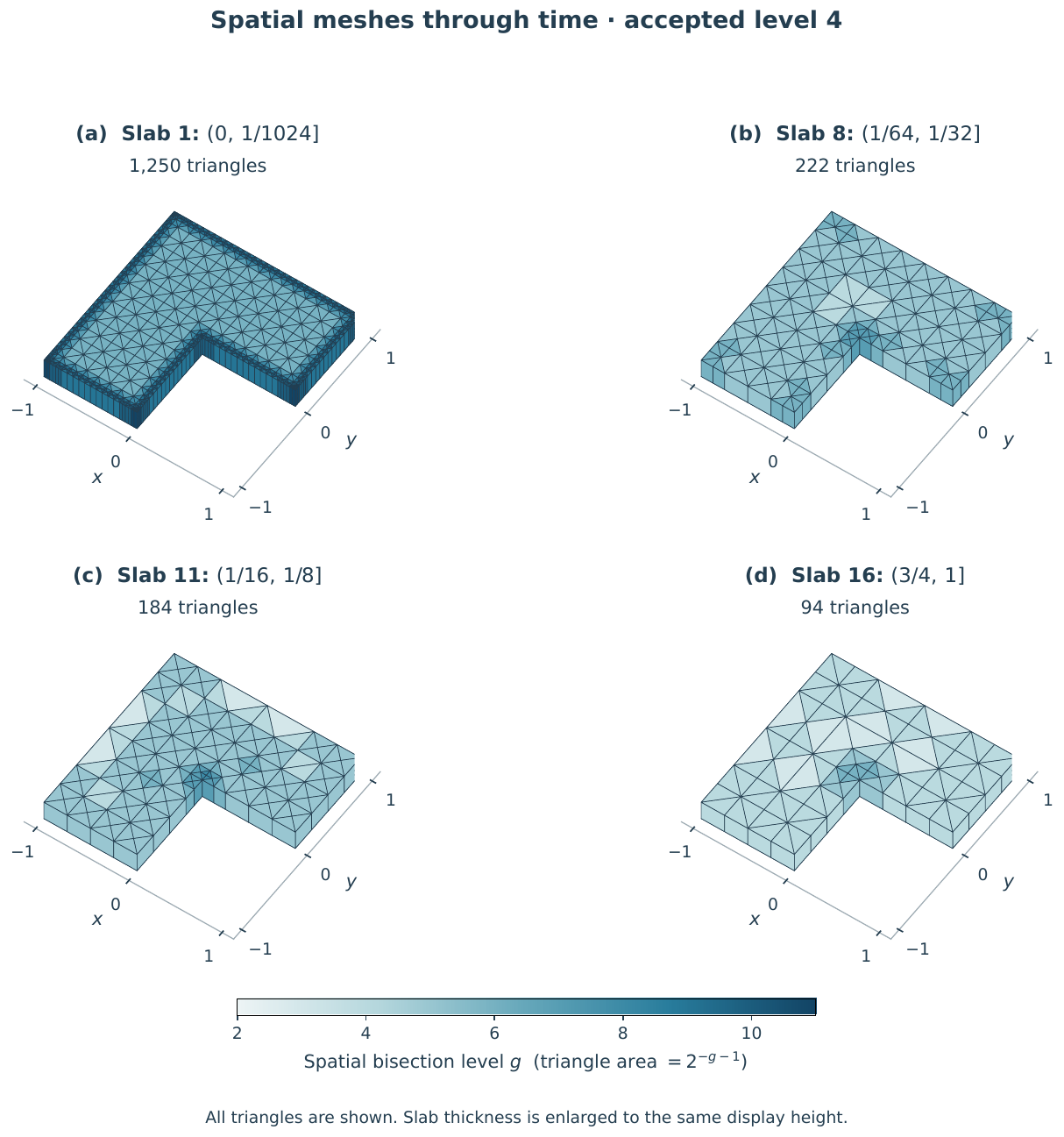}
\caption{$u_0=1$: Selected slabs of the $5764$-prism mesh: the initial layer,
two intermediate times, and the final slab. Every spatial triangle
is shown.}
\label{gen:fig-mesh-slices}
\end{figure}

Figure~\ref{gen:fig-constant} shows the estimator comparison and the
final slab profile. Figures~\ref{gen:fig-spacetime-meshes} and~\ref{gen:fig-mesh-slices}
show the actual generated meshes.

We repeat the preceding experiment with $u_0=0$, keeping everything else unchanged.

\begin{figure}
\centering
\includegraphics[width=\linewidth]{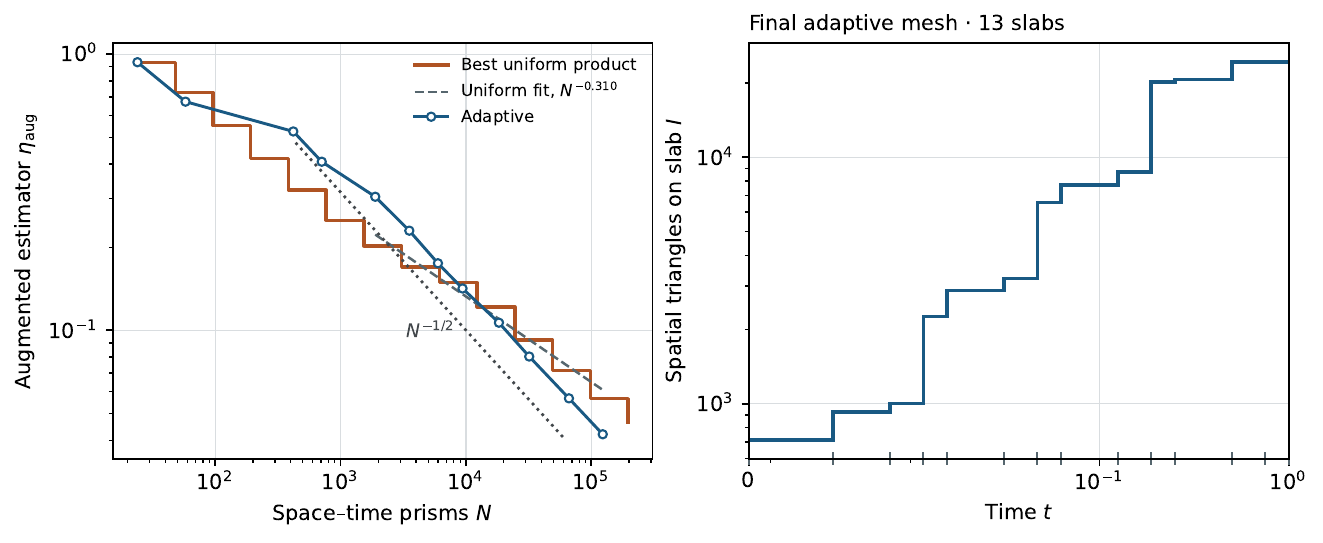}
\caption{$u_0=0$: The left panel shows
the augmented Radau estimator for the accepted meshes and the best
affordable uniform meshes.
The right panel shows the triangle count on each of the $13$ final
slabs, with marks at their endpoints.}
\label{gen:fig-zero-rates}
\end{figure}

\begin{figure}
\centering
\includegraphics[width=\linewidth,trim={0 1.7cm 0 0},clip]{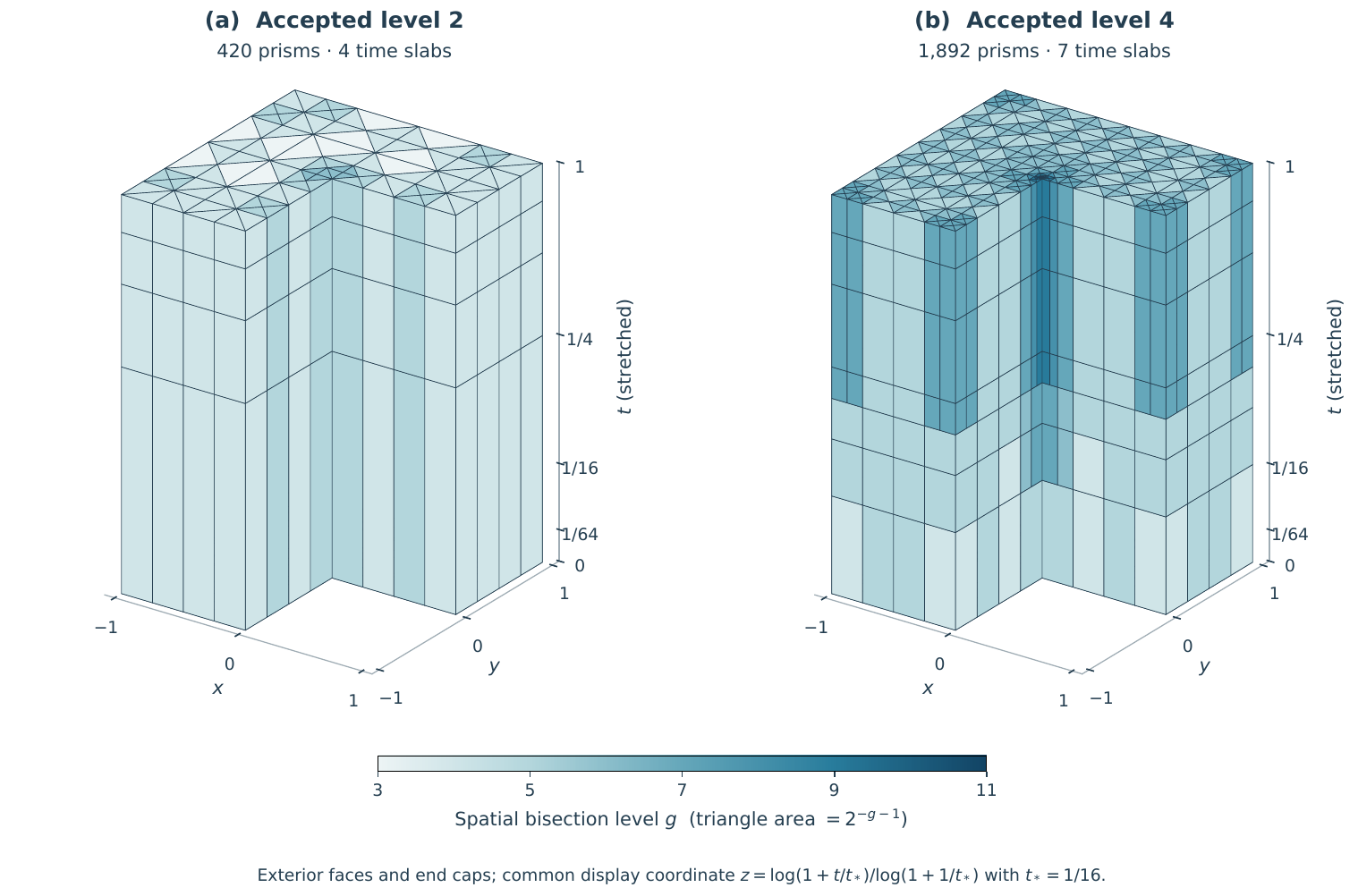}
\caption{$u_0=0$: Space--time meshes at accepted levels two and
four, with $420$ and $1892$ prisms, respectively. The color scale
shows the spatial bisection level.}
\label{gen:fig-zero-spacetime-meshes}
\end{figure}

\begin{figure}
\centering
\includegraphics[width=\linewidth,trim={0 1.7cm 0 0},clip]{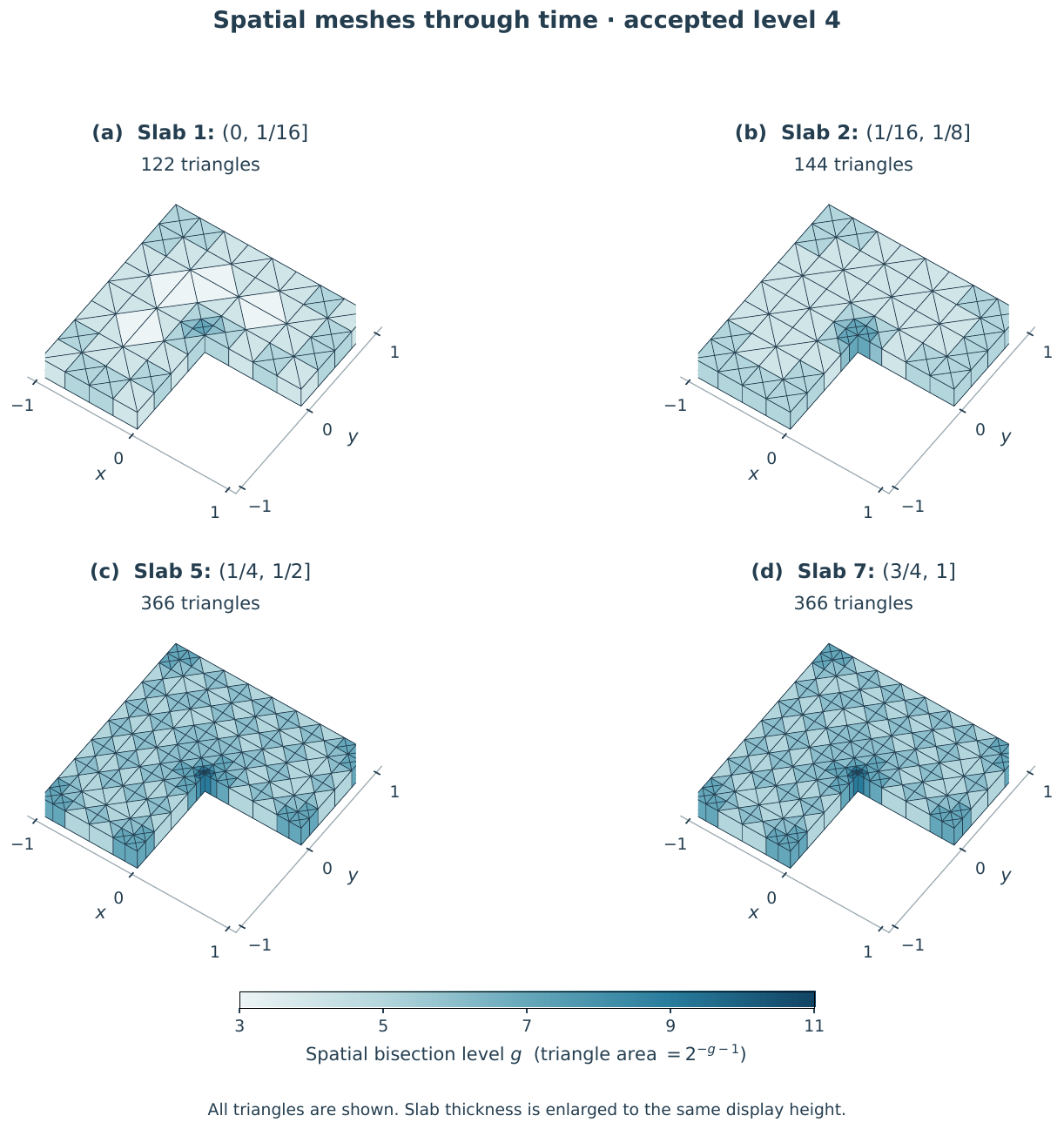}
\caption{$u_0=0$: Four slabs of the $1892$-prism mesh: the first,
second, fifth, and last. Every spatial triangle is shown.}
\label{gen:fig-zero-mesh-slices}
\end{figure}

Figure~\ref{gen:fig-zero-rates} shows the estimator comparison and
the final slab profile. Figures~\ref{gen:fig-zero-spacetime-meshes} and~\ref{gen:fig-zero-mesh-slices}
show the generated meshes at accepted levels two and four.

\bibliographystyle{alphaurl}
\bibliography{heat_tree_optimality}
\end{document}

%% file: newcommands.tex
\let\originalleft\left
\let\originalright\right
\renewcommand{\left}{\mathopen{}\mathclose\bgroup\originalleft}
\renewcommand{\right}{\aftergroup\egroup\originalright}

\newcommand{\R}{{\mathbb R}}
\newcommand{\N}{{\mathbb N}}
 
\newcommand{\Q}{\mathcal{Q}}
\newcommand{\F}{\mathcal{F}}

\newcommand{\I}{\mathcal{I}}
\newcommand{\T}{\mathcal{T}}
\newcommand{\TT}{\mathbb{T}}

\newcommand{\X}{\mathcal{X}}
\newcommand{\Y}{\mathcal{Y}}

\newcommand\SCR\ovee
\newcommand\LCC\owedge

\newcommand{\PP}[1]{\mathcal{P}^{#1}}
\renewcommand{\SS}[1]{\mathcal{S}^{#1}}

\newcommand{\QQ}{\Q}
\newcommand{\QQred}{\Q^{\LCC}}

\newcommand{\set}[2]{\left\{#1\ :\ #2\right\}}

\renewcommand{\d}{\;\mathrm{d}}

\newcommand{\Const}[1]{C_{\rm #1}}

\theoremstyle{definition}

\theoremstyle{plain}
\newtheorem{theorem}{Theorem}[section]
\newtheorem{lemma}[theorem]{Lemma}

\newtheorem{cor}[theorem]{Corollary}

\usepackage{totcount}
\newtotcounter{nconst}
\makeatletter
\newcount\rc@count
\let\rc@clearconstantlist\empty

\newcommand\rc@clearconstant[1]{\global\expandafter\let\csname rc@const@#1\endcsname\undefined}

\newcommand\resetconstants[1]{%
    \def\rc@constname{#1}
    \global\rc@count=1\relax 
    \bgroup 
        \let\\\rc@clearconstant 
        \rc@clearconstantlist
        \global\let\rc@clearconstantlist\empty 
    \egroup
}
\newcommand\const[1]{%
    \@ifundefined{rc@const@#1}{%
        \expandafter\xdef\csname rc@const@#1\endcsname{%
           \noexpand\rc@useconst{\rc@constname}{\the\rc@count}%
        }%
        \g@addto@macro\rc@clearconstantlist{\\{#1}}%
		\global\setcounter{nconst}{\the\rc@count}\relax
        \global\advance\rc@count1\relax
    }{}%
    \csname rc@const@#1\endcsname
}

\newcommand\rc@useconst[2]{\ensuremath{#1_{#2}}}
\makeatother

	\resetconstants{C}


%% file: coarsen.tex
%
%
%

\definecolor{coarsPath}{HTML}{534AB7}
\definecolor{coarsPathFill}{HTML}{CECBF6}
\definecolor{coarsCand}{HTML}{888780}
\definecolor{coarsCandFill}{HTML}{ECEAE3}
\definecolor{coarsMark}{HTML}{F5BE66}

\begin{figure}[htbp]
\centering
\footnotesize
\begin{tikzpicture}[x=10.4cm,y=0.72cm,
  cand/.style={draw=coarsCand,fill=coarsCandFill,rounded corners=1.5pt,line width=0.4pt},
  onpath/.style={draw=coarsPath,fill=coarsPathFill,rounded corners=1.5pt,line width=0.4pt},
  skipped/.style={draw=coarsCand,dashed,rounded corners=1.5pt,line width=0.5pt},
  iface/.style={draw=coarsPath,line width=1pt},
  ifacedot/.style={fill=coarsPath,circle,inner sep=1.3pt}]
  \newcommand{\timebar}[4]{%
    \path[#1] ([xshift=1pt]#3,-#2) rectangle ([xshift=-1pt]#4,-#2+0.6);}
  \node[anchor=west,inner sep=0pt] at (0,0.95)
    {\small\textbf{(a)} Cheapest path through the time tree of $\F$};
  \timebar{skipped}{0}{0}{1}
  \timebar{cand}{1}{0}{0.5}
  \timebar{onpath}{1}{0.5}{1}
  \timebar{cand}{2}{0}{0.25}
  \timebar{onpath}{2}{0.25}{0.5}
  \timebar{cand}{2}{0.5}{0.75}
  \timebar{cand}{2}{0.75}{1}
  \timebar{cand}{3}{0}{0.125}
  \timebar{onpath}{3}{0.125}{0.25}
  \timebar{cand}{3}{0.25}{0.375}
  \timebar{cand}{3}{0.375}{0.5}
  \timebar{cand}{3}{0.5}{0.625}
  \timebar{cand}{3}{0.625}{0.75}
  \timebar{onpath}{4}{0}{0.0625}
  \timebar{onpath}{4}{0.0625}{0.125}
  \timebar{cand}{4}{0.125}{0.1875}
  \timebar{cand}{4}{0.1875}{0.25}
  \node at (0.03125,-3.7) {$I_1$};
  \node at (0.09375,-3.7) {$I_2$};
  \node at (0.1875,-2.7) {$I_3$};
  \node at (0.375,-1.7) {$I_4$};
  \node at (0.75,-0.7) {$I_5$};
  \draw[iface] (0.125,-3.4) -- (0.125,-3);
  \draw[iface] (0.25,-2.4) -- (0.25,-2);
  \draw[iface] (0.5,-1.4) -- (0.5,-1);
  \node[ifacedot] at (0.0625,-3.7) {};
  \node[ifacedot] at (0.125,-3.2) {};
  \node[ifacedot] at (0.25,-2.2) {};
  \node[ifacedot] at (0.5,-1.2) {};
  \draw[line width=0.5pt] (0,-4.55) -- (1,-4.55);
  \foreach \t in {0,1}
    \draw[line width=0.5pt] (\t,-4.4) -- (\t,-4.7);
  \foreach \t in {0.0625,0.125,0.1875,0.25,0.375,0.5,0.625,0.75}
    \draw[line width=0.5pt] (\t,-4.47) -- (\t,-4.63);
  \node[anchor=north] at (0,-4.7) {$t=0$};
  \node[anchor=north] at (1,-4.7) {$t=\tend$};
  \node[anchor=north] at (0.5,-4.7) {slabs of $\F$};
  \begin{scope}[xshift=10.75cm]
    \draw[skipped] (0,0.12) rectangle +(0.33cm,0.36);
    \node[anchor=west] at (0.45cm,0.3) {skipped, exceeds $B$};
    \path[onpath] (0,-0.88) rectangle +(0.33cm,0.36);
    \node[anchor=west] at (0.45cm,-0.7) {cheapest path};
    \path[cand] (0,-1.88) rectangle +(0.33cm,0.36);
    \node[anchor=west] at (0.45cm,-1.7) {other retained intervals};
    \node[ifacedot] at (0.165cm,-2.7) {};
    \node[anchor=west] at (0.45cm,-2.7) {edge cost $K_{\F}(I,J)$};
    \node[anchor=west,inner xsep=0pt,align=left] at (0,-3.85)
      {node cost\\$L_{\F}(I)+\widehat E_I(\T_I^{\mathrm{sel}})+\beta W_I$};
  \end{scope}
\end{tikzpicture}

\vspace{1.2\baselineskip}

\begin{tikzpicture}[x=0.9cm,y=0.9cm,
  mesh/.style={draw=black!70,line width=0.4pt,line join=round},
  refedge/.style={draw=coarsPath,line width=1.2pt,line cap=round},
  chain/.style={-{Stealth[length=4pt,width=3.5pt]},line width=0.6pt},
  sw/.style={line width=0.3pt}]
  \node[anchor=west,inner sep=0pt] at (0,4.55)
    {\small\textbf{(b)} One step of the spatial search on an interval $I$};
  \fill[coarsMark] (2,2) -- (2,0) -- (1,1) -- cycle;
  \fill[coarsCandFill] (4,0) -- (2,0) -- (2,2) -- cycle;
  \fill[coarsPathFill] (4,4) -- (2,2) -- (4,0) -- (6,2) -- cycle;
  \draw[mesh] (0,0) rectangle (8,4);
  \draw[mesh] (0,0) -- (4,4) (4,0) -- (0,4)
              (2,0) -- (2,2) (0,2) -- (2,2)
              (2,0) -- (1,1) (0,2) -- (1,1)
              (4,0) -- (8,4) (8,0) -- (4,4);
  \draw[refedge] (4,0) -- (4,4);
  \draw[chain] (1.8,1) -- (2.4,1);
  \draw[chain] (2.85,0.85) -- (3.2,1.2);
  \node at (1.5,1) {$T$};
  \node at (2.62,0.42) {$T_1$};
  \node at (3.45,2.1) {$T_2$};
  \node at (4.6,2.1) {$T_3$};
  \begin{scope}[shift={(8.5,0)}]
    \node[anchor=west,inner xsep=0pt] at (0,3.75) {\textbf{Closure chain of $T$}};
    \fill[coarsMark] (0,2.95) rectangle +(0.36,0.36);
    \node[anchor=west] at (0.5,3.13) {marked leaf, maximal priority $\mathfrak p(T)$};
    \filldraw[sw,fill=coarsCandFill,draw=coarsCand] (0,2.3) rectangle +(0.36,0.36);
    \node[anchor=west] at (0.5,2.48) {coarser, must be bisected first};
    \fill[coarsPathFill] (0,1.65) rectangle +(0.36,0.36);
    \node[anchor=west] at (0.5,1.83) {executable necessary patch};
    \draw[refedge] (0,1.18) -- +(0.36,0);
    \node[anchor=west] at (0.5,1.18) {common refinement edge};
    \draw[chain] (0,0.53) -- +(0.36,0);
    \node[anchor=west] at (0.5,0.53) {crosses a refinement edge};
  \end{scope}
  \draw[chain] (4,-0.25) -- (4,-0.85);
  \node[anchor=west] at (4.15,-0.55) {bisect the patch};
  \begin{scope}[yshift=-4.59cm]
    \fill[coarsMark] (2,2) -- (2,0) -- (1,1) -- cycle;
    \fill[coarsPathFill!55] (2,0) rectangle (4,2);
    \draw[mesh] (0,0) rectangle (8,4);
    \draw[mesh] (0,0) -- (4,4) (4,0) -- (0,4)
                (2,0) -- (2,2) (0,2) -- (2,2)
                (2,0) -- (1,1) (0,2) -- (1,1)
                (4,0) -- (8,4) (8,0) -- (4,4)
                (4,0) -- (4,4);
    \draw[refedge] (2,2) -- (6,2);
    \node at (1.5,1) {$T$};
    \begin{scope}[shift={(8.5,0)}]
      \node[anchor=west,inner xsep=0pt] at (0,3.75) {\textbf{After bisecting the patch}};
      \draw[refedge] (0,3.13) -- +(0.36,0);
      \node[anchor=west] at (0.5,3.13) {new edges, $W\leftarrow W+1$};
      \fill[coarsMark] (0,2.3) rectangle +(0.36,0.36);
      \node[anchor=west] at (0.5,2.48)
        {$T$ still a leaf, $\nu(T)\leftarrow\nu(T)+1/\widehat e_I(T)$};
      \filldraw[sw,fill=coarsPathFill!55,draw=coarsPath!60] (0,1.65) rectangle +(0.36,0.36);
      \node[anchor=west] at (0.5,1.83) {next patch if $T$ is marked again};
    \end{scope}
  \end{scope}
\end{tikzpicture}
\caption{Space--time coarsening. (a)~Every interval in the temporal
refinement history of $\F$ whose spatial search stays within the budget $B$
is a node of the graph (b)~Starting from the floor mesh, the spatial
search marks a leaf $T$ of maximal priority \eqref{inst:priority}. Only
the end of the chain is bisected. The search repeats while $\max\mathfrak p>\beta$
and $W+1\le B$.}
\label{inst:fig-coarsening}
\end{figure}

%% file: heat_tree_optimality.bib
@article{inst:BD,
  author  = {Binev, P. and DeVore, R.},
  title   = {Fast computation in adaptive tree approximation},
  journal = {Numer. Math.},
  volume  = {97},
  year    = {2004},
  pages   = {193--217},
  doi     = {10.1007/s00211-003-0493-6}
}

@article{inst:BDD,
  author  = {Binev, P. and Dahmen, W. and DeVore, R.},
  title   = {Adaptive finite element methods with convergence rates},
  journal = {Numer. Math.},
  volume  = {97},
  year    = {2004},
  pages   = {219--268},
  doi     = {10.1007/s00211-003-0492-7}
}

@article{inst:BFV,
  author  = {Binev, P. and Fierro, F. and Veeser, A.},
  title   = {Near-best adaptive approximation on conforming meshes},
  journal = {Constr. Approx.},
  volume  = {57},
  year    = {2023},
  pages   = {327--349},
  doi     = {10.1007/s00365-022-09612-2}
}

@article{inst:Feischl,
  author  = {Feischl, M.},
  title   = {Inf-sup stability implies quasi-orthogonality},
  journal = {Math. Comp.},
  volume  = {91},
  year    = {2022},
  pages   = {2059--2094},
  doi     = {10.1090/mcom/3748}
}

@article{inst:Stevenson,
  author  = {Stevenson, R.},
  title   = {The completion of locally refined simplicial partitions created by bisection},
  journal = {Math. Comp.},
  volume  = {77},
  year    = {2008},
  pages   = {227--241},
  doi     = {10.1090/S0025-5718-07-01959-X}
}

@misc{FHN2026,
      title={Optimal Time-Adaptivity for Parabolic Problems}, 
      author={Michael Feischl and Fernando Henríquez and David Niederkofler},
      year={2026},
      eprint={2512.05676},
      archivePrefix={arXiv},
      primaryClass={math.NA},
      url={https://arxiv.org/abs/2512.05676}, 
}

@article{diegrading1,
  author     = {Diening, Lars and Storn, Johannes and Tscherpel, Tabea},
  title      = {On the {S}obolev and {$L^p$}-stability of the
                {$L^2$}-projection},
  journal    = {SIAM J. Numer. Anal.},
  fjournal   = {SIAM Journal on Numerical Analysis},
  volume     = {59},
  year       = {2021},
  number     = {5},
  pages      = {2571--2607},
  issn       = {0036-1429,1095-7170},
  mrclass    = {65N30 (46E35 65N12 65N50)},
  mrnumber   = {4320894},
  mrreviewer = {Nicolae\ Pop},
  doi        = {10.1137/20M1358013},
  url        = {https://doi.org/10.1137/20M1358013}
}

@misc{diegrading2,
  title         = {Grading of Triangulations Generated by Bisection},
  author        = {Lars Diening and Johannes Storn and Tabea Tscherpel},
  year          = {2023},
  eprint        = {2305.05742},
  archiveprefix = {arXiv},
  primaryclass  = {math.NA},
  url           = {https://arxiv.org/abs/2305.05742}
}

@article{actaadaptive,
  title   = {Adaptive finite element methods},
  volume  = {33},
  doi     = {10.1017/S0962492924000011},
  journal = {Acta Numerica},
  author  = {Bonito, Andrea and Canuto, Claudio and Nochetto, Ricardo
             H. and Veeser, Andreas},
  year    = {2024},
  pages   = {163–485}
}

@article{halfwavelet,
	author = {Stevenson, Rob and van Veneti{\"e}, Raymond and Westerdiep, Jan},
	date = {2022/04/06},
	doi = {10.1007/s10444-022-09930-w},
	id = {Stevenson2022},
	isbn = {1572-9044},
	journal = {Advances in Computational Mathematics},
	number = {3},
	pages = {17},
	title = {A wavelet-in-time, finite element-in-space adaptive method for parabolic evolution equations},
	url = {https://doi.org/10.1007/s10444-022-09930-w},
	volume = {48},
	year = {2022}}

@article {adaptiveG,
    AUTHOR = {Kreuzer, Christian and M\"{o}ller, Christian  and Schmidt,
              Alfred and Siebert, Kunibert },
     TITLE = {Design and convergence analysis for an adaptive discretization
              of the heat equation},
   JOURNAL = {IMA J. Numer. Anal.},
  FJOURNAL = {IMA Journal of Numerical Analysis},
    VOLUME = {32},
      YEAR = {2012},
    NUMBER = {4},
     PAGES = {1375--1403},
      ISSN = {0272-4979},
   MRCLASS = {65M60 (35K20 65M12 65M50)},
  MRNUMBER = {2991832},
MRREVIEWER = {Bruno Scheurer},
       DOI = {10.1093/imanum/drr026},
       URL = {https://doi.org/10.1093/imanum/drr026},
}

@article {adaptiveDG,
    AUTHOR = {Gaspoz, Fernando D. and Siebert, Kunibert and Kreuzer,
              Christian and Ziegler, Daniel A.},
     TITLE = {A convergent time-space adaptive {${\rm dG}(s)$} finite
              element method for parabolic problems motivated by equal error
              distribution},
   JOURNAL = {IMA J. Numer. Anal.},
  FJOURNAL = {IMA Journal of Numerical Analysis},
    VOLUME = {39},
      YEAR = {2019},
    NUMBER = {2},
     PAGES = {650--686},
      ISSN = {0272-4979},
   MRCLASS = {65M60 (65M12 65M50)},
  MRNUMBER = {3941881},
       DOI = {10.1093/imanum/dry005},
       URL = {https://doi.org/10.1093/imanum/dry005},
}

@article{FUHRER202127,
title = {Space–time least-squares finite elements for parabolic equations},
journal = {Computers and Mathematics with Applications},
volume = {92},
pages = {27-36},
year = {2021},
issn = {0898-1221},
doi = {https://doi.org/10.1016/j.camwa.2021.03.004},
url = {https://www.sciencedirect.com/science/article/pii/S0898122121000729},
author = {Thomas Führer and Michael Karkulik},
}

@article {GaSte24,
    AUTHOR = {Gantner, Gregor and Stevenson, Rob},
     TITLE = {Improved rates for a space-time {FOSLS} of parabolic {PDE}s},
   JOURNAL = {Numer. Math.},
  FJOURNAL = {Numerische Mathematik},
    VOLUME = {156},
      YEAR = {2024},
    NUMBER = {1},
     PAGES = {133--157},
      ISSN = {0029-599X,0945-3245},
   MRCLASS = {65M12 (35F35 65M15 65M50 65M60)},
  MRNUMBER = {4700409},
       DOI = {10.1007/s00211-023-01387-3},
       URL = {https://doi.org/10.1007/s00211-023-01387-3},
}

@misc{diening2025,
      title={A quasi-optimal space-time FEM with local mesh refinements for parabolic problems}, 
      author={Lars Diening and Rob Stevenson and Johannes Storn},
      year={2025},
      eprint={2502.13655},
      archivePrefix={arXiv},
      primaryClass={math.NA},
      url={https://arxiv.org/abs/2502.13655}, 
}

@article{bdd,
  author     = {Binev, P. and Dahmen, W. and DeVore, R.},
  title      = {Adaptive finite element methods with convergence rates},
  journal    = {Numer. Math.},
  fjournal   = {Numerische Mathematik},
  volume     = {97},
  year       = {2004},
  number     = {2},
  pages      = {219--268},
  issn       = {0029-599X},
  coden      = {NUMMA7},
  mrclass    = {65N50 (65N12 65N30 65Y20 68W25 68W40)},
  mrnumber   = {2050077 (2005d:65222)},
  mrreviewer = {T. Apel}
}

@article{bemao10,
  author     = {Becker, R. and Mao, S.},
  title      = {Quasi-optimality of adaptive nonconforming finite element
                methods for the {S}tokes equations},
  journal    = {SIAM J. Numer. Anal.},
  fjournal   = {SIAM Journal on Numerical Analysis},
  volume     = {49},
  year       = {2011},
  number     = {3},
  pages      = {970--991},
  issn       = {0036-1429},
  coden      = {SJNAAM},
  mrclass    = {65N30 (65N12 65N50 76D07 76M10)},
  mrnumber   = {2802555 (2012f:65185)},
  mrreviewer = {Christos A. Xenophontos}
}

@article{carstensen_2014,
  title     = {Axioms of adaptivity},
  volume    = {67},
  issn      = {0898-1221},
  url       = {http://dx.doi.org/10.1016/j.camwa.2013.12.003},
  doi       = {10.1016/j.camwa.2013.12.003},
  number    = {6},
  journal   = {Computers and Mathematics with Applications},
  publisher = {Elsevier BV},
  author    = {Carstensen, C. and Feischl, M. and Page, M. and
               Praetorius, D.},
  year      = {2014},
  month     = apr,
  pages     = {1195–1253}
}

@misc{qogeo,
      title={Quasi-orthogonality is a geometric property}, 
      author={Michael Feischl and David Niederkofler},
      year={2026},
      eprint={2611.25919},
      archivePrefix={arXiv},
      primaryClass={math.NA},
      url={https://arxiv.org/abs/2611.25919}, 
}

@article{hpopt,
author = {Canuto, Claudio and Nochetto, Ricardo H. and Stevenson, Rob and Verani, Marco},
title = {Convergence and optimality of $${\mathbf {hp}}$$hp-AFEM},
year = {2017},
issue_date = {Apr 2017},
publisher = {Springer-Verlag},
address = {Berlin, Heidelberg},
volume = {135},
number = {4},
issn = {0029-599X},
url = {https://doi.org/10.1007/s00211-016-0826-x},
doi = {10.1007/s00211-016-0826-x},
journal = {Numer. Math.},
month = apr,
pages = {1073–1119},
numpages = {47}
}

@article{ckns,
  author     = {Cascon, J. M. and Kreuzer, C. and Nochetto, R. H. and
                Siebert, K. G.},
  title      = {Quasi-optimal convergence rate for an adaptive finite
                element method},
  journal    = {SIAM J. Numer. Anal.},
  fjournal   = {SIAM Journal on Numerical Analysis},
  volume     = {46},
  year       = {2008},
  number     = {5},
  pages      = {2524--2550},
  issn       = {0036-1429},
  coden      = {SJNAAM},
  mrclass    = {65N30 (41A25)},
  mrnumber   = {2421046 (2009h:65174)},
  mrreviewer = {Hans-Peter Helfrich}
}

@article{cn,
  author  = {Cascon, J. M. and Nochetto, R. H.},
  title   = {{Quasioptimal cardinality of {AFEM} driven by nonresidual
             estimators}},
  journal = {IMA J. Numer. Anal.},
  year    = {2012},
  volume  = {32},
  number  = {1},
  pages   = {1--29}
}

@article{cr2012,
  author   = {Carstensen, C. and Rabus, H.},
  title    = {The adaptive nonconforming {FEM} for the pure displacement
              problem in linear elasticity is optimal and robust},
  journal  = {SIAM J. Numer. Anal.},
  fjournal = {SIAM Journal on Numerical Analysis},
  volume   = {50},
  year     = {2012},
  number   = {3},
  pages    = {1264--1283},
  issn     = {0036-1429},
  coden    = {SJNAAM},
  mrclass  = {74B05 (65N12 65N50 74S05)},
  mrnumber = {2970742}
}

@article{fkmp,
  author  = {Feischl, M. and Karkulik, M. and Melenk, {J.\ M.} and
             Praetorius, D.},
  title   = {Quasi-optimal convergence rate for an adaptive boundary
             element method},
  journal = {SIAM J. Numer. Anal.},
  volume  = {51},
  year    = {2013},
  number  = {},
  pages   = {1327--1348}
}

@article{gantumur,
  author  = {Tsogtgerel, G.},
  title   = {Adaptive boundary element methods with convergence rates},
  journal = { Numer. Math.},
  volume  = {124},
  year    = {2013},
  number  = {3},
  pages   = {471--516}
}

@article{infsup,
  author     = {Feischl, Michael},
  title      = {Inf-sup stability implies quasi-orthogonality},
  journal    = {Math. Comp.},
  fjournal   = {Mathematics of Computation},
  volume     = {91},
  year       = {2022},
  number     = {337},
  pages      = {2059--2094},
  issn       = {0025-5718,1088-6842},
  mrclass    = {65N30 (15A23 65N50)},
  mrnumber   = {4451456},
  mrreviewer = {Dmitrii\ Legatiuk},
  doi        = {10.1090/mcom/3748},
  url        = {https://doi.org/10.1090/mcom/3748}
}

@article{ks,
  author     = {Kreuzer, C. and Siebert, K. G.},
  title      = {Decay rates of adaptive finite elements with {D}\"orfler
                marking},
  journal    = {Numer. Math.},
  fjournal   = {Numerische Mathematik},
  volume     = {117},
  year       = {2011},
  number     = {4},
  pages      = {679--716},
  issn       = {0029-599X},
  coden      = {NUMMA7},
  mrclass    = {65N30 (65N50)},
  mrnumber   = {2776915 (2012d:65283)},
  mrreviewer = {Jeffrey S. Ovall}
}

@article{lcmhjx,
  author  = {Chen, L. and Holst, M. and Xu, J.},
  title   = {Convergence and optimality of adaptive mixed finite
             element methods},
  journal = {Math.\ Comp.},
  year    = {2009},
  number  = {265},
  volume  = {78},
  pages   = {35--53}
}

@article{ncstokes3,
  author     = {Carstensen, Carsten and Peterseim, Daniel and Rabus,
                Hella},
  title      = {Optimal adaptive nonconforming {FEM} for the {S}tokes
                problem},
  journal    = {Numer. Math.},
  fjournal   = {Numerische Mathematik},
  volume     = {123},
  year       = {2013},
  number     = {2},
  pages      = {291--308},
  issn       = {0029-599X},
  mrclass    = {65N30 (65N50 76D07 76M10)},
  mrnumber   = {3010181},
  mrreviewer = {Abdellatif Agouzal},
  doi        = {10.1007/s00211-012-0490-8},
  url        = {http://dx.doi.org/10.1007/s00211-012-0490-8}
}

@article{stevenson07,
  author     = {Stevenson, R.},
  title      = {Optimality of a standard adaptive finite element method},
  journal    = {Found. Comput. Math.},
  fjournal   = {Foundations of Computational Mathematics. The Journal of
                the Society for the Foundations of Computational
                Mathematics},
  volume     = {7},
  year       = {2007},
  number     = {2},
  pages      = {245--269},
  issn       = {1615-3375},
  mrclass    = {65N30},
  mrnumber   = {2324418 (2008i:65272)},
  mrreviewer = {Erwin Stein}
}

@article{wavelet_schwab,
  author  = {Schwab, Christoph and Stevenson, Rob},
  year    = {2009},
  month   = {09},
  pages   = {1293-1318},
  title   = {Space–time adaptive wavelet methods for parabolic
             evolution problems},
  volume  = {78},
  journal = {Math. Comput.},
  doi     = {10.1090/S0025-5718-08-02205-9}
}
